\documentclass[12pt]{article}
\usepackage{amsfonts,amssymb,amsthm,amsmath,amscd,graphicx,float}

\newcommand{\ou}{\vec{u}}
\newcommand{\ov}{\vec{v}}
\newcommand{\ow}{\vec{w}}

\newcommand{\cp}{\mathbb{C}P}
\newcommand{\cf}{\mathcal{F}_3}

\newcommand{\Fr}{Poincar\'{e} }

\newcommand{\mbQ}{\mathbb{Q}}

\newtheorem{theorem}{Theorem}
\newtheorem{lemma}{Lemma}
\newtheorem{cor}{Corollary}
\newtheorem{prop}{Proposition}
\newtheorem{defi}{Definition}

\title{Spider Evaluation and Representations of Web Groups}
\author{Charles Frohman \\ The University of Iowa}

\begin{document}

\maketitle
\begin{abstract} The topology of $SU(3)$-representation varieties of the fundamental groups of planar webs so that the meridians are sent to matrices with trace equal to $-1$ are explored, and compared to data coming from spider evaluation of the webs. Corresponding to an evaluation of a web as a spider is a rooted tree. We associate to each geodesic $\gamma$  from the root of the tree to the tip of a leaf  an irreducible component $C_{\gamma}$ of the representation variety of the web, and a graded subalgebra $A_{\gamma}$ of $H^*(C_{\gamma};\mbQ)$. The spider evaluation of geodesic $\gamma$  is the symmetrized \Fr polynomial of $A_{\gamma}$. The spider evaluation of the web is the sum of the symmetrized \Fr polynomials of the graded subalgebras associated to all maximal geodesics from the root of the tree to the leaves. \end{abstract}

\section{Introduction} A planar web $\Gamma\subset S^2$ is an oriented trivalent graph, all of whose vertices are sources or sinks, embedded in the two sphere. As $S^2\subset S^3$, it makes sense to talk about $S^3-\Gamma$. This paper is an elementary investigation of the connected components of representations of $\pi_1(S^3-\Gamma)$ into $SU(3)$. The representations we are interested in send the ``meridians'' of the web group to matrices that have trace equal to $-1$. These spaces are of interest because of their relationship to spaces used in the construction of $sl_3$-link homology, \cite{Kh1}. They are analogous to spaces of representations of knot groups into $SU(2)$ so that the meridians are sent matrices of trace $0$, that were studied in \cite{Li,RJ}.

A finite rank graded algebra $\mathcal{A}=\oplus_k A_k$ over the field $k$ is a finite direct sum of finite dimensional vector spaces, that is an associative  algebra with unit so that if $\alpha \in A_i$ and $\beta \in A_j$ then $\alpha\beta\in A_{i+j}$. The \Fr polynomial of $\mathcal{A}$ is $\sum_k\dim{A_k}t^k$. If the largest nonzero $A_k$ is $A_d$ then the symmetrized \Fr polynomial of $\mathcal{A}$ is \begin{equation}P(\mathcal{A})=t^{-d/2}\sum_k(\dim{A_k})t^k\end{equation}
We say the subalgebra $\mathcal{B}\leq \mathcal{A}$ respects the grading if
\begin{equation} \mathcal{B}=\oplus_k(\mathcal{B}\cap A_k)\end{equation} If this is the case then $\mathcal{B}$ has a symmetrized \Fr polynomial. Note the power of $t$ out in front of $P(\mathcal{B})$ might be different than $-d/2$.  Given a planar web $\Gamma$ there is a Laurent polynomial associated to $\Gamma$ via a process called spider evaluation, due to Greg Kuperberg\cite{Ku1}. The process of evaluation produces a rooted tree. A maximal geodesic in this rooted tree is an embedded path that starts with the root and ends at a leaf.  We prove that corresponding to each maximal geodesic $\gamma$  in the rooted tree there is an irreducible component $C_{\gamma}$ of the representation variety of $\Gamma$ and a subalgebra $A_{\gamma}$ of $H^*(C_{\gamma};\mathbb{Q})$ that respects the standard grading so that the sum over all geodesics of the \Fr polynomials of the $A_{\gamma}$ is equal to the spider evaluation of $\Gamma$.

Representation varieties of webs have been studied using more sophisticated tools in \cite{FKK}. That paper is a good deal more general than this one as Fontaine, Kamnitzer and Kuperberg study representations for all compact groups, though they are only concerned with the specialization of the spider polynomial at $q=1$.

A more general study of representation spaces of colored graphs was undertaken in \cite{LZ,G}.   The work of Lobb and Zentner is for $sl(N)$. However, when $N=3$ their moduli spaces coincide with the representation varieties studied here. However, we treat the irreducible components of the representation variety as separate objects, and they  study the union of all the irreducible components. This means that their computation of the \Fr polynomial has smaller coefficients than ours, as cycles in shared components only get counted once by Lobb and Zentner, and we count them with multiplicity. Spider evaluation is a sorting process, with each geodesic of the tree of resolutions picking out a different irreducible component.

The spider evaluation of the one skeleton of the cube leaves out one of the irreducible components.  Also, we  build up the component picked out by a geodesic in the resolution tree inductively starting at the leaf. Depending on the resolution, sometimes we pass to a $\cp(1)$ bundle over the component closer to the leaf,  and sometimes we are passing to a blow up.  Spider evaluation ignores the contribution to the cohomology from the exceptional divisor of the blow up. Our construction of the sublagebras of the cohomology rings of the irreducible components of the representation variety of a web just codifies these rules inductively.  There is probably a more natural construction using more sophisticated methods from algebraic geometry. However, since there are different ways of resolving a web, there is an inherently unnatural part of the process.

The goal of this paper is to give an elementary exposition of the relationship between the cohomology of the irreducible components of the representation space and a Laurent polynomial associated to the web by Kuperberg's spider evaluation.

Section 2 starts with a discussion of web groups, and the combinatorics of webs. This is followed by a description of Kuperberg's spider evaluation. The section finishes with a discussion of the geometry of the group $SU(3)$. Section 3 explores consequences of the Seifert-Van Kampen Theorem on the structure of representation varieties, and discusses blow up, and $\cp(1)$-bundles. In Section 4 we work out elementary examples of representations varieties of webs. Section 5 develops some elementary properties of representation varieties of webs that we need for our theorem. In the last section we give a method of associating to each maximal geodesic in the tree of resolvents a component of the representation variety. We show that there is a graded subring of the component associated to a geodesic whose \Fr polynomial is the summand of the spider evaluation associated to that geodesic.

 I wrote most of this paper  for graduate students, which means that in addition to being elementary, it is exploratory with many worked examples. My hope is that the paper will be readable by beginners, not just in categorification, but for people who want to study the topology of representation varieties. The author thanks Dido Salazar-Torres, Nick Teff, and Joanna Kania-Bartoszynska for extended conversations on this topic.  Finally, I thank Joel Kamnitzer for his thoughtful suggestions for how to improve the paper.

\section{Preliminaries}

\subsection{Webs, Representation Varieties, and Spider Evaluation}

A {\em web} is a planar graph with trivalent vertices and oriented
edges so that each vertex is a source or sink. 
Let $P\subset S^2\subset S^3$ be a connected planar surface and suppose that $\Gamma\subset S^2$ is a web so that $\partial P \cap \Gamma$ is contained in the edges of $\Gamma$ and the points of intersection are transverse. We call $P\cap \Gamma$ a {\em web with boundary}.  For the sake of defining the fundamental group we consider any web to be contained in $S^3$. The complement of the web is $S^3-\Gamma$.  If $\Gamma$ is a web with boundary, in the connected planar surface $P$, then $S^2-P$ is a collection of disks, $D_i$.
We consider the web to lie in the punctured $3$-sphere  $B$ which is $cl(S^3-\cup_i D_i\times [0,1])$, where the
$D_i\times [0,1]$ are suitably small collars of the $D_i$, and the $cl$ means we are taking the closure.

The fundamental group of the complement of a web is similar to a knot group. Choose a base point above the plane. For each edge of the web there
is a loop starting out at the base point, and running once around that edge so that it satisfies the right hand rule with respect to the orientation of that edge.  Such a loop is called a {\em meridian}.
If the web has $n$ edges, denote the meridians $x_i$, where $i\in \{1,2,\ldots,n\}$.  For each vertex there are three edges coming in.  There is a Wirtinger relator $x_{i_{v,1}}x_{i_{v,2}}x_{i_{v,3}}=e$ for each vertex.  The fundamental
group of the complement of the web is presented by the generators $x_i$ and the relators $x_{i_{v,1}}x_{i_{v,2}}x_{i_{v,3}}=e$.  
\begin{equation}\pi_1(S^3-\Gamma)=<x_i| x_{i_{v,q}}x_{i_{v,2}}x_{i_{v,3}}=e>\end{equation}

For instance the web shown below has three meridians, $x,y$ and $z$. Both vertices give
rise to the same relator $xyz=e$, so the group is presented as
\begin{equation} <x,y,z|xyz=e>.\end{equation}

\begin{figure}[h]\begin{center} \includegraphics{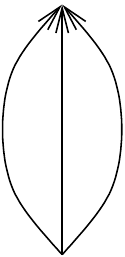}\end{center}\end{figure}

If the web has boundary, then each boundary component  of the punctured three-sphere $B$ gives rise to a punctured two-sphere in the complement of the web. 
The fundamental group of an $n$-punctured two-sphere can be presented with generators $y_i$ where $i$ ranges from $1$ to $n$ and a single relator $y_1y_2\ldots y_n=e$.  The generators are loops that encircle each puncture once counterclockwise. The inclusion map of punctured sphere into the complement of the web takes the $y_i$ to either meridians or their inverses. We orient the boundary spheres with the inward pointing normal. If an edge points inwards then the $y_i$
corresponding to the puncture gets sent to the meridian of the edge. For this reason we assign a $+$ to a boundary component of the web, so that the corresponding edge points into $B$ and a $-$ to a boundary component of the web so that the corresponding edge of the web points out of $B$. As the boundary lies on a circle we can assign a sequence of $+$'s and $-$'s corresponding to the boundary of the web, up to cyclic order.

A {\bf normal path} $\gamma$ in a web $\Gamma$ is a loop that has been pulled taut relative to its endpoints, so that it never reverses direction in an edge.  A {\bf normal loop} is a normal path where the initial and terminal edges coincide. Since the vertices of a web are all sources and sinks, as a normal loop traverses edges it alternates traveling in the direction of the edge it is in, and against the direction of the edge it is in. Thus, all normal loops have even length.

Since all the normal loops in a web $\Gamma$ have even length, the faces of a web all have an even number of sides.  A face with two sides is called a {\bf bubble}, a face with four sides is a {\bf square}. A connected web gives a cellulation of the sphere. Each edge belongs to two faces and each vertex belongs to three.  All faces have an even number of sides,
where we interpret a circle as having $0$ sides. If a face has $n$ sides, then by dividing up the edges and vertices we see it contributes
\begin{equation} \frac{6-n}{6}\end{equation} to the Euler characteristic of the sphere. Notice that if $n\geq 6$ this number is nonpositive.  Since the sum of the contributions of all faces is $2$,  we conclude  as Kuperberg did \cite{Ku2} that any connected web without boundary has a bubble or square.
Notice that the number of vertices of the graph is twice the number of edges divided by three. This means that there are an even number of vertices.  Similar computations apply to webs with boundary.

\begin{prop} If $\Gamma$ is a connected web with boundary, in a disk, and $\Gamma$ has  four boundary components, whose signs alternate as you traverse the circle  then there is an interior face of $\Gamma$ that is a square or bigon.
\end{prop} 

\qed 

Below are two graphs with four boundary components that do not have an interior square or bigon.  One of them is not connected and the signs of the boundary points of the other do not alternate. Notice there are no squares or bigons.

\begin{figure}[h]\begin{center} \includegraphics{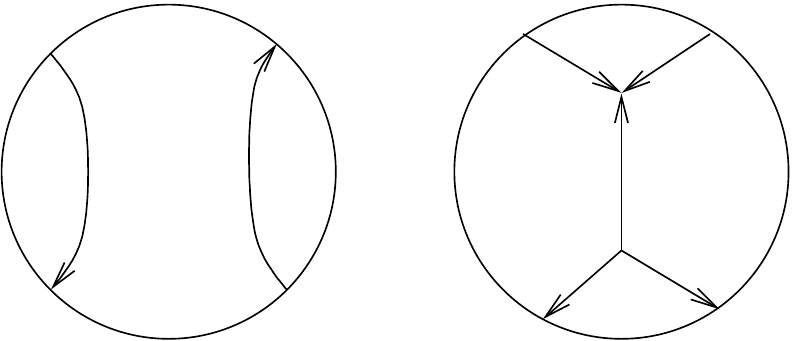} \end{center}\end{figure}

As every edge is a source or a sink, webs are {\bf bipartite graphs}.  An edge coloring of a graph is an assignment of colors to the edges so that the colors assigned to the edges at any vertex
do not repeat.  An {\bf{$r$-edge coloring}} is  an edge coloring with $r$ colors.
It is a classical theorem that if the vertices of a bipartite planar graph have maximal valence  $r$, then the graph is $r$-edge colorable \cite{M}.  That means that webs are $3$-edge colorable. 

Here is a concrete explanation of why webs are $3$-edge colorable.  Choose a face $f$ of the complement of $\Gamma$, and label it with $0\in \mathbb{Z}^3$.  If $f'$ is any other face, choose an arc $\gamma$ that starts in the interior of $f$, ends in the interior of $f'$, misses the vertices of $\Gamma$ and is transverse to the edges.
The intersection number of $\gamma$ with $\Gamma$ is the sum of the signs of the points of intersection of the oriented path $\gamma$ with the oriented edges of $\Gamma$. As $\Gamma$ defines a cycle in the $\mathbb{Z}_3$-homology of the plane, and the first homology of the plane  is zero, the intersection number is well defined as an element of $\mathbb{Z}_3$ and only depends on $f$ and $f'$.  Color the edges
between faces labeled $0$ and faces labeled $1$, by $1$. Color edges between faces labeled $1$ and faces labeled $2$ by $2$. Finally label the edges between faces labeled $2$ and faces labeled $0$ by  $0$.  Notice there are exactly three $3$-edge colorings obtainable this way and they are translates of one another.  That is, if $E$ is the set of edges of $\Gamma$ and $c,c':E\rightarrow \mathbb{Z}_3$ 
are any two colorings constructed this way $c-c'$ is a constant. We say that two colorings are equivalent if their difference is a constant. An equivalence class of colorings is a {\bf class} of colorings.
The natural class of colorings has the property that the cyclic ordering corresponding to coloring at a vertex that is a sink coincides with the counterclockwise cyclic order of the edges at that vertex.
Similarly, at a vertex that is a source the cyclic ordering from the coloring coincides with the clockwise cyclic ordering at the vertex.  Another property of the natural class of colorings is that the edges of any face are alternately colored by $2$ of the three elements of $\mathbb{Z}_3$.

\begin{prop} \label{nat} There is a natural class of $3$-edge colorings of webs by $\mathbb{Z}_3$\end{prop}

\proof The class of colorings comes from choosing a face and labeling it zero,
and then using the rule above to propagate that labeling of faces.  The edges are now labeled according the to labels of the faces to either side. Any two colorings constructed this way differ by a constant. If we reverse the orientations of all the edges simultaneously and then use the rule above to obtain a $3$-edge coloring then we get a class of colorings that are the negative of the first class.

\qed

\begin{lemma} Let $c:\Gamma\rightarrow \mathbb{Z}_3$ be a natural coloring of the connected web $\Gamma$.  For any $i\in
\mathbb{Z}_3$, let $M$ be the one manifold which is the closure of the complement of the edges $e$ with $c(e)=i$.
Then $S^2-M$ consists of a connected planar surface $P$ with $\partial P=M$ and a collection of disks. \end{lemma}

\proof A coloring $c:\Gamma\rightarrow \mathbb{Z}_3$ of a web $\Gamma$ induces an orientation at each vertex, given by the cyclic order of the edges abutting the vertex given by the coloring.  From the definition of the natural class of colorings, it is easy to see that only two colors occur on the boundary of any face.  From this you can conclude that the orientations at any two adjacent vertices of $\Gamma$ induced by the coloring are opposite to one another.  If two adjacent deleted edges lie on opposite sides of a component of $M$ the orientations of the two vertices must be the same. This means that all the deleted edges lie on the same side of any component of $M$.  

Suppose that $\Gamma$ is a connected web, and $c:\Gamma\rightarrow \mathbb{Z}_3$ is a natural coloring.  Choose a color $i \in \mathbb{Z}_3$ and consider the one manifold $M$ obtained by deleting the interiors all edges of $\Gamma$ that get assigned the the color $i$ by $c$.  This will be a collection of circles in the plane.  The circles cannot be nested, as $\Gamma$ was connected, and all the deleted edges lie on the same side of each circle.  Hence $S^2-M$ consists of a collection of open disks, and a single component that is  a planar surface $P$  whose boundary is exactly $M$.  \qed

The edges of $\Gamma$ labeled $i$ are all contained in $P$ and form a properly embedded one manifold $I\subset P$. The set $P-I$ is a collection of disks. The boundary of each disk is partitioned into two sets of  arcs and the vertices. One set of arcs is labeled $i$ and the other set of arcs is labeled by a single remaining color. Hence we can color $P-I$ with two colors. 

\begin{theorem} Up to homeomorphism connected webs correspond at most two to one with  pairs $(P,I)$ consisting of  a planar surface $P$  and a system of proper arcs $I$ embedded in $P$  so that:
\begin{itemize}
\item The union of the boundary components of $P$ and $I$ is connected.

\item The arcs $I$ separate $P$ into a collection of disks we call {\em facets}.
\item The boundary of each facet is partitioned into an even number of arcs, half of which are in $I$  and half come from arcs in the boundary of $P$.
\item The system of disks can be two colored, so that adjacent disks have the opposite color.
\end{itemize}

If two topologically distinct webs have the same underlying pair $(P,I)$ they differ only by the orientation on the edges, hence the correspondence is at most two to one.

\begin{figure}[H] \begin{center} \includegraphics{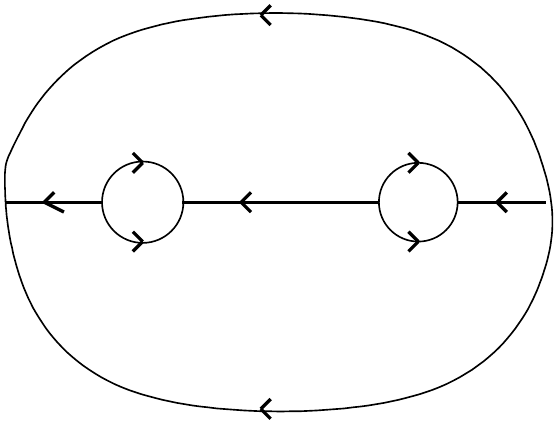} \caption{A planar web derived from a pair of pants} \end{center}\end{figure}

\end{theorem}

\qed

The three coloring corresponds to coloring the arcs in the boundary alternately, and giving the arcs in $I$ the third color.



\subsection{Spider Evaluation}

The quantized integers are Laurent polynomials in the variable $q$, 

\begin{equation} [n]=q^{-n+1}+q^{-n+3}\ldots q^{n-3}+q^{n-1}.\end{equation}

There is an invariant associated to closed webs, that takes on values in the Laurent polynomials in $q$. The invariant is  multiplicative over connected components,
based on the following rules \cite{Ku2},

\begin{itemize}
\item {\em circle}

\begin{equation} \left<\scalebox{.3}{\raisebox{-30pt}{\includegraphics{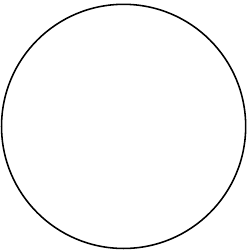}}}\right>=[3].\end{equation}

\item  {\em bubble}

\begin{equation}\label{bskein}\left<\scalebox{.3}{\raisebox{-35pt}{\includegraphics{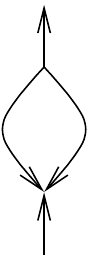}}}\right>=[2]\left<\scalebox{.3}{\raisebox{-35pt}{\includegraphics{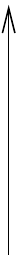}}}\right>\end{equation}

\item {\em square} 

\begin{equation}\label{sqskein}\left<\scalebox{.3}{\raisebox{-35pt}{\includegraphics{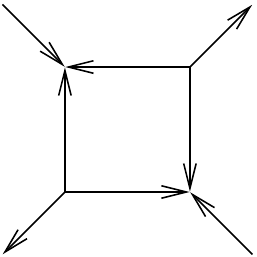}}}\right>=\left<\scalebox{.3}{\raisebox{-35pt}{\includegraphics{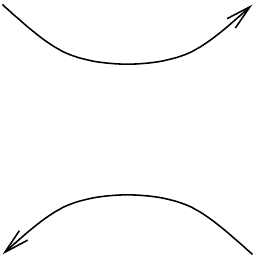}}}\right>+\left<\scalebox{.3}{\raisebox{-35pt}{\includegraphics{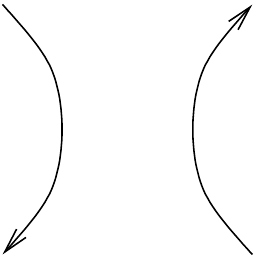}}}\right>\end{equation}
\end{itemize}

As remarked earlier, every web has a bigon or a square. Inductively collapsing when there is a bigon, or branching and smoothing when there is a square, we arrive at a collection of states with no vertices. There is a tree, with a vertex for each resolvent of the web that appears in the process of evaluation, that has a directed edge that goes from the web to the result of collapsing a bubble, or to one of the two webs obtained by deleting opposite sides of a square. We call this the {\bf tree of resolvents} $T$.  A {\bf maximal geodesic} in the tree of resolvents is a normal path that starts at the root and ends at a leaf. Since $T$ is a rooted tree, there is a unique maximal geodesic for each state. To compute the contribution of a state to the spider evaluation of $\Gamma$ start with 
$[3]^c$  where $c$ is the number of components of the state, and multiply it times $[2]^b$ where $b$ is the number of bubbles that were collapsed along the maximal geodesic corresponding to the state. The sum over all states of these contributions is the spider evaluation of the web.

For instance consider 

\begin{center}\includegraphics{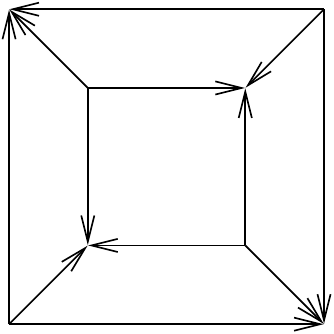}\end{center}

It has a resolution tree that looks like this.

\begin{center}\scalebox{.5}{\includegraphics{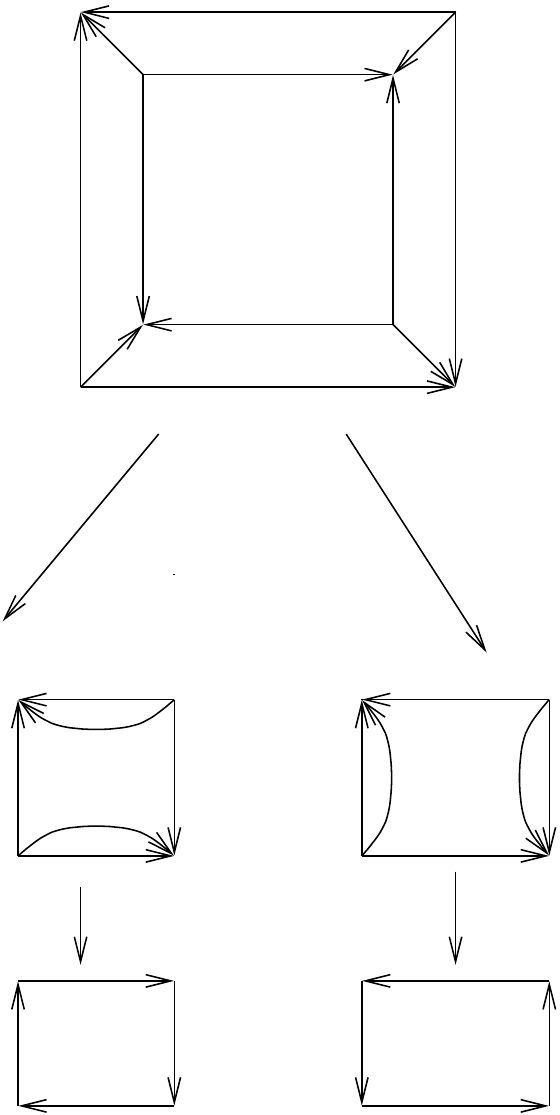}}\end{center}

The diagram is abbreviated as the end of each branch has two bubbles collapsed
at once, and technically they are collapsed one by one.    From the tree we see that the
cube evaluates as $2[2]^2[3]$.

\subsection{Getting to know $SU(3)$.}

Recall that $SU(3)$ is the group of all $3\times 3$ unitary matrices with complex coefficients having
determinant one.  The space of conjugacy classes of $SU(3)$ is a $2$-simplex. Abstractly, the maximal torus of $SU(3)$ is a two torus $S^1\times S^1$, and every element of $SU(3)$ is conjugate into that torus. However, as you can't tell what order eigenvalues come in, the conjugacy classes are in one to one correspondence with the quotient of the torus by the action of the symmetric group on three letters coming from permuting the eigenvalues.
The quotient space is a $2$-simplex.

The three vertices of the $2$-simplex
are the three central conjugacy classes, the open edges are the conjugacy classes that are homeomorphic
to $\mathbb{C}P(2)$ and the interior of the simplex consists of conjugacy classes that are homeomorphic to
$\cf$, the space of complete flags on $\mathbb{C}^3$. The stratification is given by the multiplicities of
eigenvalues.

Oddly, two elements of $SU(3)$ are conjugate if and only if they have the same trace.  This follows
from the fact that there are only two coefficients of the characteristic polynomial that count, the 
trace and the second invariant function (the sum of the pairwise products of eigenvalues).  Since the product
of all three eigenvalues is one, these two functions are complex conjugates on one another.
Hence all the the information about conjugacy is carried by one of them.

\begin{figure}\begin{center}\scalebox{.5}{\includegraphics{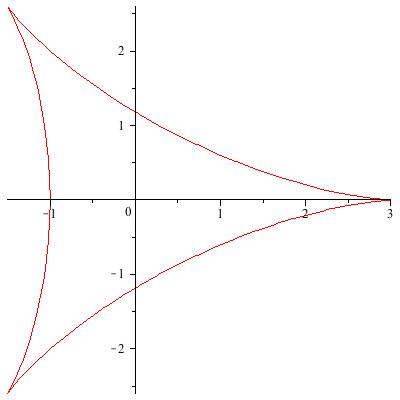}} \caption{The image of $tr:SU(3)\rightarrow \mathbb{C}$.} \end{center}\end{figure}

Notice that $SU(3)$ acts on $\mathbb{C}^3-\{\vec{0}\}$ and the action descends to $\cp(2)$.  Since
$SU(3)$ preserves the Riemannian metric on  the unit five sphere, $S^5\subset \mathbb{C}^3$
induced by the real part of the standard Hermitian pairing on $\mathbb{C}^3$, we can push it forward
to get a metric on $\cp(2)$ so that $SU(3)$ acts by isometries. This metric is called the {\em Fubini study metric}. In the Fubini Study metric, the projective plane  is a symmetric
space.  We can realize geodesic inversion $I_{<\ov>}:\cp(2)\rightarrow \cp(2)$ at $<\ov>\in \cp(2)$ as the element 
of $SU(3)$ that has
$<\ov>$ as the eigenspace of the eigenvalue $+1$ and has $-1$ as an eigenvalue of multiplicity $2$.  Since the
eigenspaces of a unitary matrix are Hermitian orthogonal to one another, this uniquely defines
the matrix.  We can identify $\cp(2)$ with this conjugacy class, which is all elements of $SU(3)$
having trace $-1$.  

It should be noted that the projective lines in the Fubini study metric are totally geodesic. Geometrically, lines 
are round spheres of radius $1/2$.  Let $<\ , \ >_H$ denote the standard Hermitian pairing.  The cosine of
the angle $\angle\ov \ow$ made by the complex lines spanned by $\ov$ and $\ow$ is,
\begin{equation} \cos{\angle \ov \ow}=\left|\frac{<\ov,\ow>_H}{||\ov||||\ow||}\right|.\end{equation}
Since the absolute value is positive we can't tell the difference between the angle $\pi$ and the angle $0$, so we just take
angles between $0$ and $\pi/2$.  This is called the {\em complex angle} between $<\ov>$ and $<\ow>$.

Complex lines $<\ov>$ and $<\ow>$ determine a unique projective line $S$ in $\cp(2)$
so long as they are not equal. The length of the short arc of the great circle spanned by $<\ov>$ and $<\ow>$ 
is  the complex angle between $<\ov>$ and $<\ow>$. If $<\ov>$ and $<\ow>$ are Hermitian orthogonal then
they are conjugate points on the projective line $S$.

The $I_{<\ov>}$ are called geodesic inversions because they fix $<\ov>$ and reverse the direction
of every geodesic running through $<\ov>$.  The transformation $I_{<\ov>}$  is the identity on a copy of $\cp(1)$ corresponding to the $-1$ eigenspace. It is also the identity map on the projective line made up of points whose underlying vectors are Hermitian orthogonal to $\ov$. 
The geodesic inversion $I_{<\ov>}$ is determined by the line $<\ov>$. It is the $+1$ eigenspace of the corresponding element of $SU(3)$, 
and $<\ov>^{\perp}$ is its $-1$ eigenspace. 

\begin{lemma}  \label{orth} The composition of two geodesic inversions $I_{<\ov>}$ and $I_{<\ow>}$ is a geodesic inversion
if and only if the vectors $\ov$ and $\ow$ are Hermitian orthogonal. \end{lemma}
\qed

The action of $SU(3)$ on $\mathbb{C}^3$ descends to an action of $SU(3)$ on $\cp(2)$. We can then look at the diagonal action of $SU(3)$ on $\cp(2)\times \cp(2)$. The orbits consist of pairs of lines
whose angle makes the same complex angle.  In general the orbits at $7$ real dimensional, except
for the two the orbit where the angle is $0$, that is the diagonal.  This means that the only invariant
sets for this action that are algebraic are the diagonal and the whole space.

\subsection{$SU(3)$-representations of web groups}

We will be considering $SU(3)$ representations of the fundamental group of the complement
of a web so that the meridians are all sent to matrices with trace equal to $-1$. 
Let $\Gamma$ be a web and let $R(\Gamma)$ be the set of homomorphisms of the fundamental group of the complement of the web into $SU(3)$ so that all the meridians are sent to matrices of trace $-1$.  Since the edges of the graph are in one to one correspondence with the meridians of the graph, the set $R(\Gamma)$ can be identified with the set of assignments of lines $<\ov>$ in $\mathbb{C}^3$ to the edges of the of the web so that the lines at any vertex are Hermitian orthogonal.  We topologize $R(\Gamma)$ as
a subspace of $\cp(2)^E$ where $E$ is the set of edges of $\Gamma$.  The spaces $R(\Gamma)$ correspond to the configuration spaces that Khovanov associates to webs in 
\cite{Kh1,KR}. Khovanov tells us that this setup originates with Kuperberg.

\begin{prop} The space $R(\Gamma)$ is in one to one correspondence to assignments of
complex lines in $\mathbb{C}^3$ to each edge of the web, so that the three lines at any vertex are Hermitian orthogonal to one another.\end{prop} \qed

\proof
 As discussed earlier such a matrix is determined by the eigenspace of the eigenvalue $1$.  By Lemma \ref{orth} if
$A,B,C\in SU(3)$, $tr(A)=tr(B)=tr(C)=-1$ and $ABC=Id$ then the $+1$-eigenspaces of $A$, $B$ and $C$ are mutually Hermitian orthogonal to one another.  The relations for the presentation
of the fundamental group of the complement of a web are all of the form $x_ix_jx_k=1$.  The matrices we are working with have order two, so $A_i$, $A_j$ and $A_k$ satisfy this relator if and only if $A_iA_j=A_k$.  This is equivalent to their $+1$-eigenspaces being Hermitian orthogonal. \qed

The space of representations of $R(\Gamma)$ of the web $\Gamma$ is at least a real  algebraic set. 
The reason I don't say complex is that the Hermitian inner product on $\mathbb{C}^3$ involves the complex conjugates of coordinate functions. However, there are homogeneous equations in the variables $z_i$ and $\overline{z}_i$ that cut out the set. {\em I have never worked an example of the representation variety of a web group, where it was obvious that some component could not be realized as a complex algebraic variety.}

 Recall a subset $V$ of a topological space $X$ is {\bf reducible} if it can be written as the union of relatively closed sets, $A$ and $B$ so that both $A-B$ and $B-A$ are nonempty. Otherwise the set is {\bf irreducible}.  A subset $V$ is an {\bf irreducible component} if it is a maximal irreducible subset of $X$.  The closure of an irreducible subset is irreducible.  An algebraic set can always be written as a finite union of irreducible components. Giving the representation variety $R(\Gamma)$ the Zariski topology from the ring of polynomials in homogeneous coordinates then it can be written uniquely as a union of irreducible components.

\section{ Fibered products, the Seifert-Van Kampen theorem and \Fr polynomials.}

Let $p:A\rightarrow X$ and $q:B\rightarrow X$ be maps.  The {\em fibered product} or fiberwise product  $A\times_{p,q}B$ of $p$ and $q$ is
  \begin{equation}(p\times q)^{-1}(\Delta)\rightarrow \Delta=X\end{equation} where 
$p\times q:A\times B\rightarrow X\times X$ is the Cartesian product of the two maps, $\Delta\subset X\times X$
is the diagonal and the map to $X$ is $p\times q$ followed by projection onto one of the factors.
The fibered product is well defined inside just about any category worth working in. For instance if the the maps are continuous, then the space $A\times_{p,q}B$ is topological and the maps are continuous.
If the maps are algebraic then the fibered product is algebraic.

\begin{prop} \label{sv} Suppose that $X\cup Y$ is a path connected topological space with $X\cap Y= A$, and that $X$ is a strong deformation retract of a path connected open set $U$, $Y$ is a strong deformation retract of a path connected open set $V$ and 
$A$ is a strong deformation retract of the path connected open set $U\cap V$. Let  $p:R(X)\rightarrow R(A)$ and
$q:R(Y)\rightarrow R(A)$ be the restriction maps on the space of representations of the fundamental groups of these spaces into some group $G$.
The restriction map \begin{equation}R(X\cup_A Y)\rightarrow R(A)\end{equation} is the fibered product of $p$ and $q$.\end{prop}

\proof  Recall that $\rho_1:\pi_1(X)\rightarrow G$ and $\rho_2:\pi(Y)\rightarrow G$ correspond to
a representation $\rho:\pi_1(X\cup Y)\rightarrow G$ if and only if $\rho_1|_{\pi_1(A)}=\rho_2|_{\pi_1(A)}$.
Hence the space of homomorphisms of $\pi_1(X\cup Y)$ into $G$ can be identified with the set of
ordered pairs of homomorphism $(\rho_1,\rho_2)$ whose restriction to $\pi_1(A)$ agrees.  This is
the definition of the fibered product. \qed

\begin{prop}If $p:A\rightarrow X$ is a fiber bundle and $q:A\rightarrow X$ is also a fiber bundle  then the  fibered product of $A$ and $B$ is a fiber bundle with fiber equal to the Cartesian product
of the fibers of $p$ and $q$.\end{prop}

When we are working examples we will encounter maps that are  blowups. 
An exposition of the blowup of a complex
manifold along a complex submanifold can be found in \cite{GH}.

If $Y\subset X$ is a complex codimension $2$ submanifold then at each point of $Y$ there is a coordinate
patch $(z_1,\ldots,z_n)$ so that in these coordinates $Y$ is given by the equations $z_1=z_2=0$.  Consider
the set

\begin{equation}B= \{(z_1,z_2,\ldots,z_n,[x_1,x_2])\in \mathbb{C}^n\times \cp(1)|z_1x_2=z_2x_1\}\end{equation}

First note the set is well defined because the equation is homogeneous in $x_1$ and $x_2$.  

The projection map $\pi:B\rightarrow \mathbb{C}^n$ ,  given by \begin{equation}\pi(z_1,\ldots,z_n,[x_1,x_2])=(z_1,\ldots,z_n)\end{equation} has the property that the inverse image of any point where $z_1\neq 0$ or $z_2\neq 0$ is a single point, except when $z_1=z_2=0$ the inverse image is a copy of $\cp(1)$.  Locally $B$ is the blowup of $\mathbb{C}^n$ along $z_1=z_2=0$.  The copy of $\cp(1)$ over any point $y\in Y$ can be identified with the projective line modeled on the normal space to $T_yY$ in $T_yX$.

The blow up of a complex $n$-manifold is locally a complex  $n$-dimensional submanifold of $\mathbb{C}^n\times \cp(1)$ as
the equations are nondegenerate.

To define the blow up along a codimension $2$ complex submanifold $Y$ of the complex manifold $X$, cover $Y$ with open neighborhoods in $X$ that are coordinate charts like above, construct it locally and glue it together. The submanifold $Y$ is the {\em singular locus} of $p$ and its inverse image is the {\em exceptional divisor}

We are specifically interested in the blow up of $\cp(2)\times \cp(2)$ along the diagonal $\Delta$.
Consider $\tilde{\cp(2)\times\cp(2)}\subset \cp(2)\times \cp(2)\times \cp(2)$ which is the set of all triples $([\ou],[\ov],[\ow])$
so that $\ow$ is Hermitian orthogonal to $\ou$ and $\ov$. This can be seen to be an algebraic set, by replacing $\ow$ by its complex conjugate and requiring the conjugate to be complex orthogonal to the first two lines.  The projection map $p:\tilde{\cp(2)\times \cp(2)}\rightarrow \cp(2)\times \cp(2)$ can be seen to be equivalent to the blowup of $\cp(2)\times \cp(2)$ along the diagonal.

Suppose that $p:\tilde{X}\rightarrow X$ is a blow up along the complex codimension $2$ submanifold $Y$, and $E=p^{-1}(Y)$ is the exceptional divisor. The cohomology $H^*(\tilde{X};\mbQ)$ is a vector space over the ring $H^*(X;\mbQ)$ via 
\begin{equation} p^*:H^*(X,\mbQ)\rightarrow H^*(\tilde{X};\mbQ).\end{equation} The vector space has an easy to understand structure.  The submanifold $E$ defines a line bundle $[E]$ over $\tilde{X}$ This line bundle has a first Chern class $c([E])\in H^2(\tilde{X};\mbQ)$. The pullback
$p^*H^*(X;\mbQ)\rightarrow H^*(\tilde{X};\mbQ)$ is injective, and
\begin{equation} H^*(\tilde{X};\mbQ)=p^*H^*(X;\mbQ)\oplus c([E])\cup p^*H^*(Y;\mbQ).\end{equation}
Furthermore \begin{equation}c([E])\cup p^*H^*(Y;\mbQ)\cong H^*(E;\mbQ)/p^*H^*(Y;\mbQ)=H^{*-2}(Y;\mbQ).\end{equation}
As the \Fr polynomial only sees the free part,  we can use the rational cohomology groups to
compute Betti numbers.
This justifies the formula for symmetrized \Fr polynomials,
\begin{equation} P(\tilde{X})=P(X)+P(Y).\end{equation}
Notice that $p^*H^*(X;\mbQ)\leq H^*(\tilde{X})$ is an isomorphic copy of the homology of $X$.

We won't always be working with blowups along a submanifold, but we will come very close.  Let
\begin{equation}\tilde{\cp(2)\times\cp(2)}=\{([\ou],[\ov],[\ow])\in \cp(2)^3|\ou,\ov,\ow \ \text{are mutually Hermitian orthogonal}\}.\end{equation} The projection map $p:\tilde{\cp(2)\times\cp(2)}\rightarrow \cp(2)\times \cp(2)$ defined
by $p([\ou],[\ov],[\ow])=([\ou],[\ov])$ is topologically equivalent to the blowup of $\cp(2)\times \cp(2)$
along its diagonal $\Delta=\{([\ou],[\ou])\}$.  We denote $p^{-1}(\Delta)$ by $\tilde{\Delta}$. It is an oriented $\cp(1)$ bundle over $\Delta$.  Let $X$ be an algebraic set. Suppose that $q:X\rightarrow \cp(2)\times \cp(2)$ is onto and algebraic. The fibered product $\tilde{X}$ of $p$ and $q$ acts very much like a blowup. We call the projection map $b:\tilde{X}\rightarrow X$ a {\bf generalized blow up}.
\begin{equation} \tilde{X}=\{(x,y)\in X \times \tilde{\cp(2)\times \cp(2)}|q(x)=p(y)\},\end{equation}
with projection map $b:\tilde{X}\rightarrow X$ given by $b((x,y)=x$.  It should be noted that
$b$ is algebraic and onto. The fibers of $b$ away from the {\em singular locus} $Y=q^{-1}(\Delta)$ consist of a single point,
and the restriction of $b$ to $\tilde{Y}=b^{-1}(Y)$ is an oriented $\cp(1)$-bundle over $Y$.  It is still the case that $p^*:H^*(X;\mbQ)\rightarrow H^*(\tilde{X};\mbQ)$ is an injective ring homomorphism. It should be noted that if $q:X\rightarrow \cp(2)\times \cp(2)$ is algebraic then the generalized blow up $\tilde{X}$ is algebraic. 

\begin{theorem} If $b:\tilde{X}\rightarrow X$ is projection from the  fibered product of $p:\tilde{\cp(2)\times\cp(2)}\rightarrow \cp(2)\times \cp(2)$ with $q:X\rightarrow \cp(2)\times \cp(2)$ where $q$ is an onto algebraic mapping then $b_*:H_*(\tilde{X})\rightarrow H_*(X)$ is onto. \end{theorem}

\proof  Let $\Delta\subset \cp(2)\times \cp(2)$ be the diagonal. Let $S\subset \tilde{X}$ be the inverse image of $\Delta$ under the projection map of the fibered product of $p$ and $q$. Since $S$ is at least a closed real algebraic subvariety of  an algebraic variety it is an NDR. That means there exists $\tilde{U}$ an open neighborhood of $S$ so that $S$ is a strong deformation retract of $\tilde{U}$. Since $b:\tilde{X}-S\rightarrow X-q^{-1}(\Delta)$ is a homeomorphism $U=b(\tilde{U})$ is open, and the deformation retraction of $\tilde{U}$ onto $S$ descends to give a deformation retraction of $U$ onto $q^{-1}(\Delta)$.

We work with the Mayer-Vietoris sequence for $\tilde{X}$ based on $\tilde{X}-S,\tilde{U}$, and for $X$ based on $X-q^{-1}(\Delta),U$. The map $b$ induces a map between the two sequences.

Notice that $b:\tilde{U}-S\rightarrow U_q^{-1}(\Delta)$ and $b:\tilde{X}-S\rightarrow X-q^{-1}(\Delta)$ are  homeomorphisms so the induced maps on homology are isomorphism.  The map $b:S\rightarrow q^{-1}(\Delta)$ is an oriented $\cp(1)$-bundle so that it is onto , so the map induced by $b:\tilde{U}\rightarrow U$ on homology is onto. The rest is a diagram chasing argument. \qed

\begin{cor} If $p:\tilde{X}\rightarrow X$ is a generalized blow up, then  $p^*:H^*(X;\mbQ)\rightarrow H^*(\tilde{X};\mbQ)$ is an injective ring homomorphism.\end{cor}

\proof Since for any space $B$, $H^*(B;\mbQ)=Lin(H_*(B);\mbQ)$ we have that $p^*$ is the adjoint of an onto mapping. It is therefore injective. \qed

\begin{prop} \label{blowex} The fibered product of a blowup of a codimension 2 submanifold of a manifold with itself
has two irreducible components. A $\cp(1)\times\cp(1)$ bundle over the singular locus and
a copy of the original blow up.\end{prop}

\proof Suppose now that $\pi:W\rightarrow X$ is the blowup along $Y$.  Let's form the fibered product with itself.

This means first form $\pi\times \pi :W\times W\rightarrow X\times X$ and let $FP=(\pi\times\pi)^{-1}(\Delta)$.
The projection map to $\Delta=X$ where $\Delta$ is the diagonal of $X\times X$ is just the restriction of the projection map $\pi \times \pi$.  We have 
$p:FP\rightarrow X$ which is the fibered product of  the blow up with itself.  Notice that $p^{-1}(Y)$ is a complex $n$-dimensional manifold that is a fiber bundle over $Y$ with fiber $\cp(1)\times \cp(1)$.  Notice that $p^{-1}(X-Y)$ is also a complex $n$-dimensional manifold with fiber a point.  In local coordinates,
\begin{equation}p^{-1}(z_1,\ldots,z_n)=(z_1,\ldots,z_n,[x_1,x_2],z_1,\ldots,z_n,[x_1',x_2'])\end{equation}
where $z_1x_2=z_2x_1$ and $z_1x_2'=z_2x_1'$, and not both $z_1$ and $z_2$ are zero. In $\cp(1)$ there
is only one solution to this equation so in fact

\begin{equation}p^{-1}(z_1,\ldots,z_n)=(z_1,\ldots,z_n,[x_1,x_2],z_1,\ldots,z_n,[x_1,x_2])\end{equation}

Notice if  a sequence of points like this is converging to something over $Y$, in local coordinates the sequence
converges to a point in the diagonal of $\cp(1)\times\cp(1)$.  The closure of $p^{-1}(X-Y)$ does not
contain $p^{-1}(Y)$ so the fibered product of the blowup with itself is reducible. However $p^{-1}(Y)$
and the closure of $p^{-1}(X-Y)$ are manifolds so they are irreducible. \qed

A similarly easy situation is when $p:E\rightarrow B$ is an {\bf orientable } $\cp(1)$-bundle.
Orientable means that you can pick an atlas of local trivializations where the change of coordinates
induce orientation preserving homeomorphisms of $\cp(1)$.  This is always true of $\cp(1)$-bundles that occur in an algebraic setting. The Leray-Hirsch theorem \cite{HA} says that if $p:E\rightarrow B$ is an orientable $\cp(1)$-bundle over a paracompact base $B$ then there is a class $\gamma \in H^2(E)$ whose restriction to each fiber is the first Chern class of the canonical line bundle over the fiber. The projection map $p:E\rightarrow B$ makes $H^*(E;\mbQ)$ a vector space on basis $1,\gamma$ over $H^*(B;\mbQ)$.
 Hence  $P(E)=[2]P(B)$.  We summarize:

\begin{theorem} If $p:E\rightarrow B$ is the projection map of either an oriented $\cp(1)$-bundle then $p^*:H^*(B;\mbQ)\rightarrow H^*(E;\mbQ)$ is an injective ring homomorphism. \end{theorem}

\section{Elementary Examples of Representation Varieties of Webs}

We begin with studying the representations varieties of the simplest webs.

\subsection{A Circle} Suppose that $\Gamma$ is a single circle.  From the definition, we see that
$R(\Gamma)$ is just $\cp(2)$.  The circle can be oriented counterclockwise or clockwise. The symmetrized \Fr polynomial of $\cp(2)$ is $[3]$.

\vspace{.1in}

If a web $\Gamma$ has boundary, it lies inside of a planar surface $P$. We can orient the boundary
components of $P$ and then assign a cyclic sequence of signs to each boundary component by traversing it in a positive direction and noting whether the edge of $\Gamma$ points in or out.
if the edge points into $P$ its a $+$ and if it points out its a $-$.

\subsection{A Triad}

There are actually two triads, one that is a source and one that is a sink.  A triad has three boundary components that are $(+,+,+)$ or $(-,-,-)$ as you traverse the boundary of the disk.  I have pictured the $(+,+,+)$ version.

\begin{center}\includegraphics{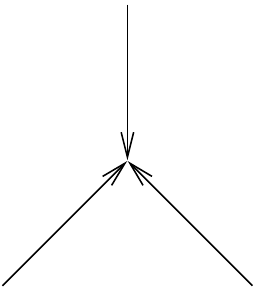}\end{center}

Its representation variety is the space of complete flags $\cf$ on $\mathbb{C}^3$.  To make this identification,
choose a base vertex and a sense. The first line is the one dimensional subspace of the flag, the span of the first two is the two dimensional subspace of the flag. The projection of the representation variety into any two boundary values determines it, as there is a unique line orthogonal to any plane in $\mathbb{C}^3$.  The symmetrized \Fr polynomial of $\cf$ is $[2][3]$.

\subsection{A Bubble}\label{bubs}

\begin{center}\includegraphics{bubble.pdf}\end{center}

The bubble is a $\cp(1)$ bundle over $\cp(2)$ sitting as the diagonal in $\cp(2)\times \cp(2)$ which
you can think of as the representations of the boundary (the diagonal that is). The representation space is just $\cf$.

Here is how to see it. We are assigning lines to each edge of the bubble, so that all the lines at each vertex are Hermitian orthogonal.  If we assign $[\vec{u}]$ to the bottom edge, and $[\vec{v}]$ and $[\vec{w}]$ to the sides of the bubble, then $<\vec{v},\vec{w}>$ span the perpendicular to $[\vec{u}]$.
As we are only in three dimensional space, the requirement that the line assigned to the top edge be orthogonal to $<\vec{v},\vec{w}>$ means that it is $[\vec{u}]$.  Hence the representation variety is determined by the lines assigned to the bottom edge, and the left edge of the bubble, $([\vec{u}],[\vec{v}])$ where $\vec{u}$ and $\vec{v}$ are Hermitian orthogonal.  The projection map to the value on the bottom edge is a fiber bundle, with fiber the copy of $\cp(1)$ that is perpendicular to $[\vec{u}]$ in $\cp(2)$.

\subsection{The Jumping Jack}

\begin{center}\includegraphics{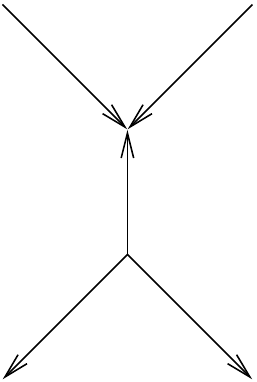}\end{center}

The jumping jack has four boundary components and two vertices.  The four boundary components are arranged $(+,+,-,-)$ as you go around the circle.  The $+$ means it points into the disk and the $-$ means the edge points out.  The representation variety has a single component which is a $\cp(1)$-bundle over $\mathbb{F}_3$. If you
project it into $\cp(2)^2$ corresponding to picking  the values on two edges going to the two $+$ boundary components
or the two $-$ boundary components the image is a copy of $\mathbb{F}_3$ embedded in $\cp(2)$ by sending a flag
to a line spanning the  line in the flag, followed by a line that is Hermitian orthogonal to the first line so the two span the
two dimensional subspace in the flag.  The map is a fiber bundle with fiber $\cp(1)$. However,
it should be noted that this last mapping is not holomorphic.  If instead you project onto the values
of two adjacent edges of opposite sign, then the map is a blow-down.  In this case make sure that the edges correspond to the first edge at each vertex so that the map is holomorphic. The singular locus is the diagonal in $\cp(2)^2$.
For this reason we denote the space of representations of the jumping jack by $\widetilde{\cp(2)\times \cp(2)}$.

\subsection{A Square}

\begin{center}\includegraphics{square.pdf}\end{center}

The boundary of the square has four points arranged $(+,-,+,-)$ as you go around the circle.
The representation variety of the square has two components.  The image of the two components under restriction to the boundary yields two different sets \begin{equation}\Delta_v=\{(<\ou>,<\ou>,<\ov>,<\ov>)\in \cp(2)^4|\ov,\ou \in \mathbb{C}^3\}\end{equation}
and \begin{equation}\Delta_h=\{(<\ou>,<\ov>,<\ov>,<\ou>)\in \cp(2)^4|\ov,\ou \in \mathbb{C}^3\}\end{equation}  sitting in $\cp(2)^4$.  
We refer to the sets $\Delta_h$ and $\Delta_v$ as the components of the {\em large diagonal} of $\cp(2)^4$. If we start at the lower left hand corner and proceed counterclockwise to enumerate the boundary components, then the lines assigned to the vertices in $\Delta_v$ agree on vertices that are vertically aligned, and the lines assigned to the vertices in $\Delta_h$ agree if they are horizontally aligned.

The
set \begin{equation}\delta=\{(<\ou>,<\ou>,<\ou>,<\ou>)\in \cp(2)^4|\ou \in \mathbb{C}^3\}\end{equation} will be called
the {\em small diagonal}.
Each component is the blow up of $\Delta_v$ or $\Delta_h$ along the small diagonal.
Hence each component is a copy of $\widetilde{\cp(2)\times \cp(2)}$.  The intersection of the two components is a copy of $\mathcal{F}_3$.  

\begin{defi} \label{square}Whenever we have an evaluation map to the corners of the square we use $\Delta_h$, $\Delta_v$ and $\delta$ to denote the sets as above. \end{defi}

In the figure below I have indicated the values in $\cp(2)^4$  corresponding to $\Delta_v$ on the outside, and the assignments corresponding to $\Delta_h$ on the inside.

\begin{center}\begin{picture}(130,130) \includegraphics{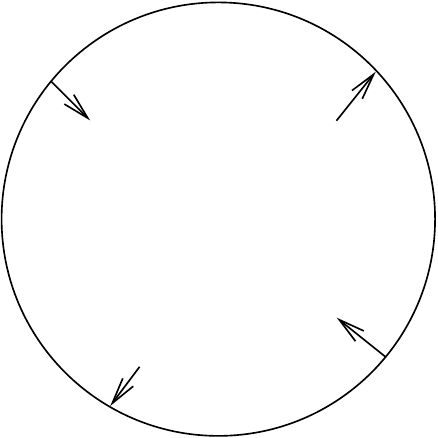}\put(-10,10){$\ou$}
\put(-25,35){$\ou$} \put(-10,110){$\ou$}\put(-25,85){$\ov$} \put(-120,110){$\ov$}\put(-110,85){$\ov$} \put(-120,10){$\ov$}\put(-110,35){$\ou$}\end{picture}\end{center}

Here is a careful analysis of the representations of the square.

\begin{center}\begin{picture}(148,148)  \scalebox{2}{\includegraphics{square.pdf}} \put(0,0){$<\ou>$}\put(-20,74){$<\ov>$}\put(-95,20){$<\ow>$}\put(-110,128){$<\alpha\ou+\beta\ow>$}
\put(-180,74){$<\gamma\ou+\delta\ov>$}\end{picture}\end{center}

We place an orthogonal frame at the lower right vertex. This means that the upper horizontal edge  has a line
assigned to it, that is Hermitian orthogonal to $\ov$. Following
the diagram around the bottom, the line assigned to the left hand vertical side must be orthogonal to $\ow$ so it gets a value in the space of $\ou$ and $\ov$.  Now its crunch time as the upper horizontal edge and the left vertical edge meet at the upper left vertex.  From this we see that either $\gamma=0$ or $\alpha=0$. This gives rise to two irreducible components of
the space of representations of the complement of the square.

\begin{itemize}
\item The left and right vertical edges are assigned the same line,
\item The top and bottom vertical edges are assigned the same line.
\end{itemize}

In the first case  we are free to choose any vector that is orthogonal to $\ow$ and the rest is determined as shown.

\begin{center}\begin{picture}(148,148)  \scalebox{2}{\includegraphics{square.pdf}} \put(0,0){$<\ou>$}\put(-20,74){$<\ov>$}\put(0,148){$<\ou>$}\put(-95,20){$<\ow>$}\put(-95,128){$<\ow>$}
\put(-160,74){$<\ov'>$}\put(-182,0){$<\ou '>$} \put(-182,148){$<\ou'>$}\end{picture}\end{center}

This means that the space is a $\cp(1)$-bundle (from the choice of $\ou'$) over $\mathcal{F}_3$ ( from the lower right corner).

The other component is similar except that the top corners , and lower corners agree.

There is another way of viewing the components of the representation space of the square.
Suppose you know the values on the edges coming in from the corners. From our analysis
of the representation space, either the corners that are over each other, or the corners that
are across from each other are equal.  These two situations correspond to the image of the component lying in $\Delta_v$ or $\Delta_h$ under restriction to the boundary values. Suppose for the moment that the two pairs of edges over one another are equal.  Suppose the pair on the left are assigned $[\ou']$ and the pair on the right are assigned $<\ou>$. If $<\ou>\neq [\ou']$ then the lines assigned horizontal sides of the square
are determined by the fact that they must be perpendicular to both $<\ou>$ and $[\ou']$.  
The vertical sides are now forced because they share a vertex with the horizontal sides and the edges from the corners. However if $<\ou>=[\ou']$ then there is a $\cp(1)$'s worth of
choices for the lines assigned to the horizontal sides as they must only be orthogonal to 
to $<\ou>$.  The lines assigned to the vertical sides are still determined. From this we see
that the representation space is the result of blowing up $\cp(2)\times \cp(2)$ along
the diagonal.  Notice that the closure of $\cp(2)\times \cp(2)-\Delta$ is dense in this space.

\subsection{The Double Square} The representation variety of a double square has some interesting features.

\begin{center}\includegraphics{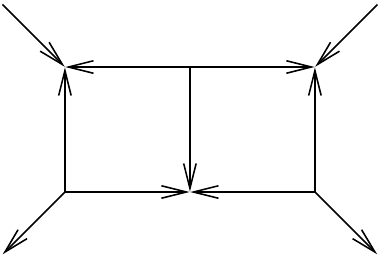}\end{center}

The most expedient way of working this example is to see it as a jumping jack on the left glued to a square on the right, over the two horizontal edges on the left.  Let $Ev_D:R(doublesquare)\rightarrow \cp(2)^2$ be the evaluation map onto these two edges.  The representation variety of the jumping jack has one component which is the blow up of $\cp(2)^2$ along the diagonal under $Ev_D$.  The representation variety of the square has two components. One is
a blow up over the diagonal under $Ev_D$ and the other is an iterated $\cp(1)$-bundle over the diagonal of $\cp(2)^2$ under $Ev_D$.  As the representation variety of the double square is the fibered product of these, we get three components. The first is a blow up of the diagonal
of $\cp(2)$ and the other two are iterated $\cp(1)$ bundles over the diagonal of $\cp(2)$.
These two components are interesting as one is vertical over the left square, but horizontal over the right square, and the other is vertical over the right square and horizontal over the left square. They are diffeomorphic via a canonical choice of map, as each space can be identified with a component of the representation variety of the square obtained by smoothing.

Let $\theta:R(doublesquare)\rightarrow R(jumping \ jack)$ the the map that forgets the assignment
of lines to the two vertical edges that form the sides of the right square. Notice that the component
of the representation of the double square that is the blow up of the diagonal, maps homeomorphically to the representation variety of the jumping jack.  This phenomenon occurs whenever you have a family of squares lying side by side.  The component of the representation variety where all the vertical edges carry the same label, maps homeomorphically down the representation variety of the web obtained by deleting any two of the vertical edges of the squares.

\subsection{The Cube}

An important example to understand is the representation variety of the one skeleton of the cube.

\begin{center}\includegraphics{cube.pdf}\end{center}

It has a nonseparating square, which you can think of as the inner square in the diagram. If you excise its four sides, you get another square, which is the outer square in the diagram. Choose the first
edge at each vertex to be the diagonal edge. That way the evaluation map to the four edges,
\begin{equation} Ev_D:R(cube)\rightarrow \cp(2)^4\end{equation} is holomorphic. By the Seifert-Van Kampen theorem the representation
variety of the cube is the fibered product of the representation variety of the square with itself. The representation variety of the square has two components, corresponding to the blow ups of the two components of the large diagonal $\Delta_v$ and $\Delta_h$ along the small diagonal.  There are three components of the representation variety of the cube. Two components are copies of the components of the representation variety of the square, one mapping onto $\Delta_v$ and one mapping onto $\Delta_h$ as blowups. The third component is vertical. Its is a $\cp(1)\times \cp(1)$ bundle over the small diagonal.
The third component appears four times. Twice from taking the fibered product of a blow up with itself, and twice from taking the fibered product of components of the representation variety with different images under evaluation.

Here is the interesting part about the cube.  All three components are diffeomorphic.
The vertical component can be seen by labeling the square above so that all three diagonal edges
get the same line $<\ou>$. Now label the edges of the outer square alternately $<\ov>$, $<\ow>$ where $\ou,\ov,\ow$ is an Hermitian orthonormal frame.  Similarly label the edges of the inner square
$\ov',\ow'$ where $\ou,\ov',\ow'$ is an Hermitian orthonormal frame. Projecting to $\cp(2)^4$ by taking the values on the diagonal edges, this maps onto the small diagonal and the fiber is $\cp(1)\times \cp(1)$ coming from the choice of two Hermitian orthonormal frames including $\ou$.  However you can cut the cube into two squares by taking a vertical line or a horizontal line.  Notice that this component under those two projections to $\cp(2)^4$ is a blow up of a component of the large diagonal over the small diagonal.  In fact, for all three components, there is a way of dividing the cube into two squares so that the component projects down to the small diagonal.   Recalling the spider evaluation of the cube, it is the sum of the symmetrized Poincar\'{e} polynomials of two of the components.  This means that the mode of evaluation ignores one of the components.  This is a feature of polynomial invariants of knots, links and graphs associated to representations of quantum groups.  Somehow the planarity of the diagram is needed to sort out the data that counts, and the data that doesn't.

\subsection{The Double Hexagon}

\begin{center} \includegraphics{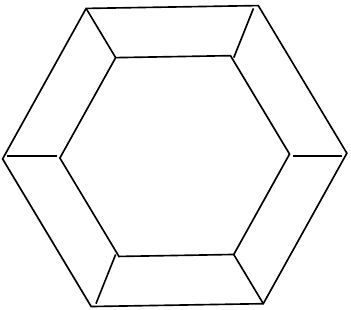} \end{center}

All of the components of the representation varieties that we have explored up to now are manifolds.  This is only the case for a small class of examples.
The representation variety of the double hexagon is the simplest example I know of a representation variety with a singular component. There is a component of its representation variety which is the result of blowing up $\cp(2)^3$ three times along subvarieties corresponding to points where the first two coordinates are the same, the second two coordinates are the same and the first and third coordinates are the same. The component of the representaition variety has a singularity over the locus in $\cp(2)^3$ consisting of points where all  three coordinates are the same. 

\section{Some Observations on the Structure of the Representation Variety of a Web}

\subsection{Forgetful Maps}

\begin{figure}\begin{center}\includegraphics{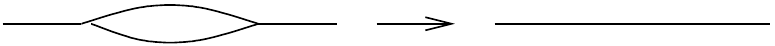} \end{center} \caption{Smashing a bubble}\label{smashb} \end{figure}

In the process of resolving a web into states, we make two kinds of moves.
The first is smashing a bubble shown in Figure \ref{smashb}.

\begin{figure}\begin{center}\includegraphics{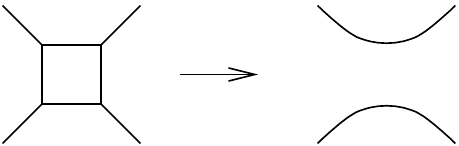}\end{center}\caption{Smoothing a square}\label{sqs}\end{figure}

The second is smoothing a square, shown in Figure \ref{sqs}. There are actually two ways to smooth a square, we just show the horizontal smoothing.

Given an irreducible component of the result of one of these moves, we need to find an irreducible component of its anticedent that gets mapped to it. To do this we need a clear picture of how representations of the two webs correspond. We call the maps we use to do this  {\bf forgetful maps}.  Recall that if $\Gamma$ is a web with edges $E$, then $R(\Gamma)$ can be identified with functions $\rho:E\rightarrow \cp(2)$, so that if two edges $e,e'$ share a vertex $\rho(e)$ is Hermitian orthogonal to $\rho(e')$.  

\begin{itemize}
\item Suppose that $\Gamma'$ is obtained from $\Gamma$ by smashing a bubble.
  Let $E$ be the set of edges of $\Gamma$ and let $e_0,e_1,e_2,e_3$ the the edges of $\Gamma$ that touch the bubble. Here we imagine $e_0$ going into the bubble, $e_1,e_2$ the sides of the bubble and $e_3$ exiting the bubble. Let $e_4$ be the edge of $\Gamma'$ that replaces the four edges of the bubble. If $E'$ is the edges of $\Gamma'$ we have
  \begin{equation} E'=E-\{e_0,e_1,e_2,e_3\}\cup\{e_4\}.\end{equation}
  Any representation of $\Gamma$ restricts to a representation of the bubble. From our analysis of the representations of a bubble if $\rho:E\rightarrow \cp(2)$ is a representation, $\rho(e_0)=\rho(e_3)$. Define a map $f:E'\rightarrow E$ by letting $f$ take any edge the two webs share in common to itself, and sending $e_4$ to $e_0$.  The forgetful map
  \begin{equation} \phi:R(\Gamma)\rightarrow R(\Gamma')\end{equation} is given by $\phi(\rho)=\rho\circ f$.

\item Suppose now that $\Gamma'$ is obtained from $\Gamma$ by smoothing a square. Let $E$ be the edges of $\Gamma$ and suppose $e_{tl},e_{tm},e_{tr},e_{vl},e_{vr},e_{bl},e_{bm}$ and $e_{br}$ are the edges that touch the squares.  The indices are suggestive of their positions. For instance $e_{tl}$ the edge coming in from the top left, and $e_{br}$ is the edge leaving on the bottom right. Let $e_{t}$ and $e_{b}$ be the edges of $\Gamma'$ the replace the square. Here you should think of $e_t$ replacing the top three sides of the square, $e_{b}$ replacing the bottom three sides of the square, and the two vertical sides of the square are deleted. If $E'$ is the edges of $\Gamma'$ then
  \begin{equation} E'=E-\{e_{tl},e_{tm},e_{tr},e_{vl},e_{vr},e_{bl},e_{bm},e_{br}\}\cup \{e_t,e_b\}.\end{equation}

  Any representation of $\Gamma$ restricts to a representation of the square. From our analysis of representations of the square, if $\rho\in R(\Gamma)$ either the pairs $\{e_{tl},e_{tr}\}$ and $\{e_{bl},e_{br}\}$ take on the same values or
  the pairs  $\{e_{tl},e_{bl}\}$ and $\{e_{tr},e_{br}\}$ take on the same values.
  Let $R(\Gamma)'$ be the subset of the representation variety where  $\{e_{tl},e_{tr}\}$ and $\{e_{bl},e_{br}\}$ take on the same values.

  Let $f:E'\rightarrow E$ be the map that takes edges that the two webs share in common to themselves and sends $e_t$ to $e_{tl}$ and $e_b$ to $e_{bl}$. In this case the forgetful map
  \begin{equation} \phi:R(\Gamma)'\rightarrow R(\Gamma)\end{equation} is given by $\phi(\rho)=\rho\circ f$.

  \end{itemize}
 
\subsection{Foundations}

\begin{lemma} A closed web has nonempty representation space and any edge takes on all values \end{lemma}
\proof The proof is by induction on the number of vertices. From Euler characteristic, a web always has a square or bubble.  If there is a bubble then the representation variety of the web is a $\cp(1)$-bundle over the representation variety of the web with the bubble smashed. By the inductive hypothesis this is nonempty.  

Suppose that $\Gamma$ has a square. We can assume that $\Gamma$ is connected otherwise we can work component by component.  Let $\Gamma_h$ and $\Gamma_v$ be the two webs obtained from $\Gamma$ by smoothing. At least one of $\Gamma_h$ and $\Gamma_v$ is connected. Without loss of generality assume it $\Gamma_h$.  Notice $\Gamma_h$ is obtained from $\Gamma$ by deleting two edges.  A representation of $R(\Gamma_h)$ allows us to assign lines to every
edge of $\Gamma$ except for the four sides of the square. By the inductive hypothesis $R(\Gamma_h)$ is nonempty. Let $[\ou]$ and $[\ov]$ be the lines assigned to the two horizontal edges in $R(\Gamma_h)$.  We use these to assign vectors to the diagonals going into the square. There is always $[\ow]$ that is Hermitian orthogonal to both $[\ou]$ and $[\ov]$, assign this to the two vertical sides of the square.  Finally, there is a unique way to fill out the horizontal sides of the square so that all lines at each vertex of the square are Hermitian orthogonal to one another. 

Now that we know that $R(\Gamma)$ is nonempty, the action on of $SU(3)$ on
$R(\Gamma)$ by conjugation looks like the action of $SU(3)$ by translation on
$\cp(2)$ when you restrict to the values on any edge, so every edge takes on all values.  \qed

\begin{lemma} \label{twob} Let $\Gamma$ be a web with two boundary components with signs $(+,-)$. Let $Ev_B:R(\Gamma)\rightarrow \cp(2)^2$ be the map that takes any representation to the two lines assigned to the edges coincident with the boundary.  $Ev_B(R(\Gamma))$ is the diagonal $\Delta$ of $\cp(2)^2$ and
\begin{equation} Ev_B:R(\Gamma)\rightarrow \Delta\end{equation} is a fiber bundle. Furthermore if $\Gamma'$ is the closed web
obtained adding an edge to the boundary  of $\Gamma$ then $R(\Gamma)=R(\Gamma')$.\end{lemma}

\proof The proof is by  induction on the number of vertices. The inductive step comes from comparing with the result of collapsing a bubble or smoothing a  square. The base case is a single line segment. 
Clearly the evaluation map takes this to the diagonal.  Once again, Euler characteristic guarantees there is always a square or bubble, and the representations of the web $\Gamma$ correspond to representations of the simpler web in a way that does not effect boundary values.  
Don't forget that our representations are into  $SU(3)$ and they are closed under conjugation. By conjugating by elements of $SU(3)$ we see we get the full diagonal, and the fibers over any two points are the same and the map is locally trivializable.

There is a one to one correspondence between representations of $\Gamma$ and $\Gamma'$ as adding the edge to the outside does not restrict the values that can be assigned to the two boundary edges. \qed

\begin{prop} \label{fourb} Let $\Gamma$ be a web with four boundary components having signs $(+,-,+,-)$ as you traverse the circle and let $Ev_D:R(\Gamma)\rightarrow \cp(2)^4$ be the evaluation map.
If $C$ is any component of $R(\Gamma)$ then $Ev_D(C)\subset \Delta_v$ or $Ev_D(C)\subset \Delta_h$.
If  $Ev_D(C)\subset \Delta_v$ then $C$ is diffeomorphic
to a component of $R(\Gamma_v)$, the closed web, obtained by closing the outside
as below on the left. The correspondence inducing the diffeomorphism just takes assignments of edges to assignments of edges.  If  $Ev_D(C)\subset \Delta_h$ then $C$ is diffeomorphic to a component of $R(\Gamma_v)$, the result of closing up $\Gamma$ as shown on the right below. \end{prop}

\proof This is  similar to the last proof. The base case is that the web consists of two line segments. The representation variety is $\cp(2)\times \cp(2)$ so it only has one irreducible component $C$. Depending on how we order the vertices $Ev_D(C)\subset \Delta_h$ or $Ev_D(C)\subset \Delta_v$.

Suppose Proposition  has been proved for all webs with four boundary components with signs  $(+,-,+,-)$ and less than $n$ vertices, and $\Gamma$ has $n>0$ vertices, and $C$ is a component of $R(\Gamma)$. If $\Gamma$ is not connected the result is true as a consequence of Proposition \ref{twob}. If $\Gamma$ is connected, from Euler characteristic considerations $\Gamma$ has a square or a bigon.  There is a web $\Gamma'$ obtained by smashing the bigon or smoothing the square appropriately, and an irreducible component of $R(\Gamma')$ so that the forgetful map $q:R(\Gamma)'\rightarrow R(\Gamma')$ is onto.  By induction $Ev_D(C')\subset \Delta_v$ or  $Ev_D(C')\subset \Delta_h$. The map $Ev_D:C\rightarrow \cp(2)^4$ factors through $Ev_D:C'\rightarrow \cp(2)^4$ via the forgetful map. Therefore the proposition is true for $C$.   \qed

\begin{center}\includegraphics{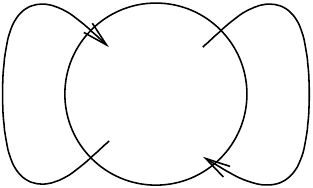} \hspace {1in} \raisebox{-.25in}{\includegraphics{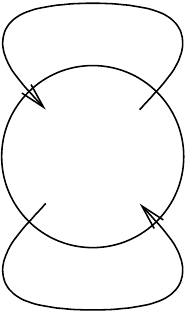}}\end{center}

\subsection{All components of the representation variety are big or small at each square} A pair of edges $\bf{e}=\{e_1,e_2\}$ of the web $\Gamma$ forms and {\bf{isthmus}} if there is a disk $D\subset S^2$ that $\partial D$ intersects $\Gamma$ transversely, and there are exactly two points of intersection, one in the interior of $e_1$ and one in the interior of $e_2$. This separates $\Gamma$ into two webs with boundary, each having two boundary components.  From the lemma above, and  Proposition \ref{sv}  we conclude that the evaluation map
\begin{equation} Ev_{\bf{e}}:R(\Gamma)\rightarrow \cp(2)^2\end{equation}
that takes each representation to its values on $e_1$ and $e_2$ has image the diagonal of $\cp(2)^2$.

We call a square $S$  in a web $\Gamma$ {\bf separating} if two of the edges adjoining $S$ form an isthmus.  Necessarily these two edges abut $S$ along one side. In the process of spider evaluation, applying the relation from Equation \ref{sqskein}, the result of deleting the two sides of the square adjacent to the side of the square that touches the isthmus creates a web with one more boundary component than you started with which we denote $\Gamma'$ and one with the same number of components $\Gamma''$.

\begin{defi}\label{hor} Let $\Gamma$ be a web and $S\subset \Gamma$ a square.  Let
\begin{equation} Ev_S:R(\Gamma)\rightarrow \cp(2)^4\end{equation} be the map that takes a representation to lines assigned
to the four diagonals of the square.
 There is  map $R(\Gamma)\rightarrow R(\Gamma')$ that comes from forgetting the lines assigned to the four deleted edges.  Let $C'$ be the image of $C$ under this map.
If the image of $C'$ is one of the components of the large diagonal of $\cp(2)^4$ then we say $C$ is {\bf big}.
If the image of $C'$ is the small diagonal we call $C'$ {\bf small}.
\end{defi}

Notice that all the components of the representation variety of the square are big.  For the cube, two components are big and one is small with respect to any face of the cube, though the same component can be big with respect to one face and small with respect to another.

\begin{prop}\label{biggy}  Let $\Gamma$ be a web with four boundary components having signs $(+,-,+,-)$ as you traverse the circle and let $Ev_D:R(\Gamma)\rightarrow \cp(2)^4$ be the evaluation map.  Let $Ev_D:R(\Gamma)\rightarrow \cp(2)^4$ be the evaluation map at the boundary. Let $C$ be an irreducible component of $R(\Gamma)$ then $C$ is  big  in the sense that $Ev_D:C\rightarrow \cp(2)^4$ either has image $\Delta_h$, or $\Delta_v$ . \end{prop}

\proof The proof will be by induction on the number of vertices of $\Gamma$. In the base case, $\Gamma$ is just two arcs and the Proposition is true, by the analysis at the beginning of the proof of  Proposition \ref{fourb}. Suppose now that the Proposition has been proved true for all webs with less than $n$ vertices, and suppose that $\Gamma$ is a web with $n$ vertices, and $C$ is an irreducible component of $R(\Gamma)$. If $\Gamma$ is not connected the conclusion of the proposition holds by Proposition \ref{twob}.

If $\Gamma$ is connected then it has a bigon or a square.  If $\Gamma$ has a bigon, let $\Gamma'$ be the result of smashing the bigon and $\phi :R(\Gamma)\rightarrow R(\Gamma')$ be the forgetful map. The component $C$ is a $\cp(1)$ bundle over an irreducible component $C'$ of $R(\Gamma')$. The evaluation map $Ev_D:C\rightarrow \cp(2)^4$ factors through $Ev_D:C'\rightarrow \cp(2)^4$ via $\phi$. Since $\Gamma'$ has fewer vertices than $\Gamma$ the inductive hypothesis allows us to conclude that $Ev_D:C'\rightarrow \cp(2)$ is one of the large diagonals or the small diagonal.  Therefore the conclusion of the hypothesis holds for $C$.

If $\Gamma$ has a square $S$, place a disk $D'$ about $S$ and let $D''$ be the disk with the same boundary as $D'$ so that $D'\cap D''=S^2$,  the sphere that $\Gamma$ lies in.  By Proposition \ref{fourb} there is a way of adding two arcs to $\Gamma$ to get a closed web $\Gamma'$ so that $C$ can be identified with an irreducible component of $R(\Gamma')$.  Let $\Gamma''=\Gamma'\cap D''$, Since $\Gamma''$ has fewer vertices than $\Gamma$, the conclusion of the proposition holds for $\Gamma''$ with respect to $Ev_S:C\rightarrow \cp(2)^4$, where $Ev_S$ is evaluation at the four corners of the square.  This means there is a smoothing of the square $S$, $\Gamma'''$ so that the forgetful map $\phi:R(\Gamma)'''\rightarrow R(\Gamma''')$ takes $C$ onto an irreducible component $C'''$ of $R(\Gamma''')$. Since $\Gamma'''$ has fewer vertices than $\Gamma$ the conclusion of the proposition holds for $\Gamma'''$, the image of $\Gamma'''$ is either $\Delta_h$, $\Delta_v$ or $\delta$. Since $Ev_D:C\rightarrow \cp(2)^4$ factors through $Ev_D:C'''\rightarrow \cp(2)^4$ by the forgetful map, we have that the conclusions of the proposition hold for $C$. \qed

On the other hand

\begin{prop} Let $\Gamma$ be a connected web, and suppose that $\Gamma$ has a square $S$. 

\begin{itemize}
\item  If $S$ is separating then
for one of the big diagonals, $\Delta_v$ or $\Delta_h$ there exists
a component $C$ of $R(\Gamma)$ that is big and onto that component.  The other component
of the large diagonal
has the property that it's intersection with the image of $Ev_S$ is the small diagonal.

\item If $S$ is nonseparating then for each component $\Delta_h$ or $\Delta_v$ of the large diagonal
of $\cp(2)^4$ there is at least one component $C$ of $R(\Gamma)$ that is big with respect to the square whose image is that component. 

\end{itemize}
\end{prop}

\proof  Excising the four sides of the square, we get a web $\Gamma'$ with four boundary components
in a disk, that are labeled $(+,-,+,-)$ as you traverse its boundary.  There is a map $R(\Gamma)\rightarrow R(\Gamma')$ coming from forgetting the values of the representation on the four sides of the square.  The map is onto and takes components to components, so in analyzing the values of $Ev_S$ it makes no difference if we study $R(\Gamma)$ or $R(\Gamma')$.  

If $S$ is separating then $\Gamma'$ has two components. The representation variety of $\Gamma'$ is the Cartesian product of the representation variety of its two components, so we get the desired result by applying Proposition \ref{twob}.

If $S$ is nonseparating, by Euler characteristic considerations $\Gamma'$ has a square or bigon in the
interior of the disk, using Proposition \ref{fourb} we see that each component of the representation
variety of $\Gamma'$ maps onto a component of the web with the bigon collapsed, or onto one of
the smoothings of the web. Hence we can proceed by induction on the number of vertices to get the desired result.

\qed

Putting things together, given a square $S$ that is nonseparating, then
for each smoothing $\Gamma_h$ and $\Gamma_v$ there is a one to one correspondence between components of $R(\Gamma)$ that are big with respect to $S$ and
have image $\Delta_h$ (respectively $\Delta_v$) and components of $R(\Gamma_h)$ (respectively $R(\Gamma_v)$).  

If the square $S$ is separating suppose that
$\Gamma_h$ has one more component than $\Gamma$ and $\Gamma_v$ has the same number of components as $\Gamma$. There is a one to one correspondence between components of $R(\Gamma)$ that  are big at $S$ and components of $R(\Gamma_h)$.  The mapping induced by restriction from the big component of $R(\Gamma)$  $C$ is mapped to the component $C'$ of $R(\Gamma_h)$ is a generalized blow up along the intersection of $C'$ with the small diagonal.  The small components of $R(\Gamma)$ are orientable $\cp(1)$ bundles over $R(\Gamma_v)$.

\section{The Correspondence between Spider Evaluation and the Cohomology of the Representation Variety}

Let $\Gamma$ be a web. Corresponding to a maximal geodesic is a sequence of webs,
\begin{equation} \Gamma=\Gamma_0,\Gamma_1,\ldots,\Gamma_n\end{equation}
where each web $\Gamma_i$ is obtained from $\Gamma_{i-1}$ by either collapsing a bubble, or deleting two sides of a square. 

We now describe an inductive procedure for choosing an irreducible component $C_i$ of the representation varieties of the $\Gamma_i$ along with a graded subalgebra $A_i\leq H^*(C_i;\mbQ)$.

If $\Gamma_n$ is a leaf, then it consists of a disjoint union of circles. The representation variety has a single irreducible component that is a Cartesian product of copies of $\cp(2)$, one for each circle. Hence $C_n=R(\Gamma_n;\mbQ)$ and $A_n=H^*(R(\Gamma_n))$.

\begin{itemize}
\item If $\Gamma_i$ is obtained from $\Gamma_{i-1}$ by a bubble move, and we have chosen a component $C_i$ of the representation variety of $\Gamma_i$ along with a graded subalgebra $A_i$ of $H^*(\Gamma_i;\mbQ)$, there is a unique component $C_{i-1}$ of the representation variety of $\Gamma_{i-1}$ that is an oriented $\cp(1)$-bundle over $C_i$ Let $\phi_i:C_i\rightarrow C_{i-1}$ be the restriction of the forgetful map.  If $\gamma\in H^2(C_{i-1})$ restricts to the first Chern class of the canonical bundle of each fiber. Define
  \begin{equation} A_i=\phi_i^*(A_{i-1})<1,\gamma>\leq H^*(A_i)\end{equation}
  be the subalgebra of $H^*(A_i;\mbQ)$ generated by $\phi_i^*(A_{i-1})$ and $\gamma$.

\item Suppose that  $\Gamma_i$ is obtained from $\Gamma_{i-1}$  by smoothing a square $S$. Let $C_i$ be the component of $R(\Gamma_i)$ that has been constructed. Let $D$ be a disk containing the square $S$ so that $\partial D$ intersects the four edges coming out of $S$, and let $D'$ be the disk with the same boundary as $D$ so that $D\cup D'=S^2$ the plane containing $\Gamma_{i-1}$. Let $\Gamma_{i-1}'=\Gamma_{i-1}\cap D'=\Gamma_i\cap D'$. This last identification is because the part where  $\Gamma_i$ and $\Gamma_{i-1}$ differ is contained in $D$.

  We can identify $C_i$ with a component of $R(\Gamma_{i-1}')$ by Proposition \ref{fourb}.   By Proposition \ref{biggy}, $Ev_S:C_i\rightarrow \cp(2)^4$ is either a big diagonal or the small diagonal.

  \begin{enumerate}
  \item In the case that the image $Ev_S:C_i\rightarrow \cp(2)^4$ onto a large diagonal, so that $C_i$ is big , there is a component $C_{i-1}$, so that  the forgetful map restricted to $C_{i-1}$,  $\phi_{i}:C_{i-1}\rightarrow C_i$ is the generalized  blow up along the small diagonal. Let $A_{i-1}=\phi_i^*(A_{i})$.
  \item In the case that $Ev_S:C_i\rightarrow \cp(2)^4$ is the small diagonal, there is a unique component $C_{i-1}$ of $R(\Gamma_{i-1})$ so that the forgetful map restricted to $C_{i-1}$,  $\phi_i:C_{i-1}\rightarrow C_i$ is an oriented $\cp(1)$-bundle. Let $A_{i-1}=\phi_i^*(A_i)$.
\end{enumerate}

\end{itemize}

\begin{theorem} Let  $\Gamma$ be a planar web. Let $T$ be a tree of resolvents. For each maximal geodesic $\gamma$ in $T$ let $C_{\gamma}$ be the irreducible component of $R(\Gamma)$ picked out by this procedure. Let $A_{\gamma}$ be the graded subalgebra of $H^*(C_{\gamma})$ constructed by the procedure above. The symmetrized \Fr polynomial of $A_{\gamma}$ is the contribution of that geodesic to the spider evaluation of $\Gamma$. Therefore the spider evaluation of $\Gamma$ is the direct sum of the symmetrized \Fr polynomials of the $A_{\gamma}$. \end{theorem}
\qed

\proof Notice that the \Fr polynomial of the algebra $A_{\gamma}$ of any geodesic is $[2]^b[3]^c$ where $b$ is the number of bubbles that were smashed along the geodesic and $c$ is the number of components of the state at the end of the geodesic. This coincides with its contribution to the spider evaluation of $\Gamma$.  \qed

\end{document}